\documentclass[12pt]{amsart}
\usepackage{amssymb}
\usepackage[dvips]{graphics}
\ifx\mathrm\undefined\let\mathrm\rm\fi
\ifx\mathbf\undefined\let\mathbf\bf\fi
\ifx\mathfrak\undefined\let\mathfrak\frak\fi
\ifx\mathcal\undefined\let\mathcal\cal\fi
\ifx\mathbb\undefined\let\mathbb\Bbb\fi
\ifx\emph\undefined\let\emph\it\fi
 at9.98pt

\newcommand{\g}{{{\mathfrak g}\,}}

\newcommand{\n}{{{\mathfrak n}}}

\newcommand{\h}{{{\mathfrak h\,}}}

\newcommand{\R}{{\mathbb R}}
\newcommand{\C}{{\mathbb C}}

\newcommand{\la}{\lambda}

\newcommand{\dontprint}[1]
{\relax}

\newtheorem%
{thm}{Theorem}[section]
\newtheorem%
{proposition}[thm]{Proposition}
\newtheorem%
{lemma}[thm]{Lemma}
\newtheorem%
{lemmadef}[thm]{Lemma-Definition}
\newtheorem%
{corollary}[thm]{Corollary}
\newtheorem%
{conjecture}[thm]{Conjecture}

\newcommand{\bea}{\begin{eqnarray*}}
\newcommand{\eea}{\end{eqnarray*}}
\newcommand{\bean}{\begin{eqnarray}}
\newcommand{\eean}{\end{eqnarray}}

\newcommand{\nc}{\newcommand}
\nc{\on}{\operatorname}
\nc{\al}{\alpha}
\nc{\ri}{\rangle}
\nc{\lef}{\langle}
\nc{\W}{{\mathcal W}}
\nc{\La}{\Lambda}
\nc{\ep}{\epsilon}
\nc{\Om}{\Omega}
\newcommand{\be}{\begin{displaymath}}
\newcommand{\ee}{\end{displaymath}}
\newcommand{\bs}{\boldsymbol}
\nc{\PCr}{{ \Bbb P  (\C[x])^r   }}
\newtheorem{theorem}{Theorem}[section]

\newcommand{\id}{{{\rm id}}}
\newcommand{\A}{{{\mathcal{C}}}}
\newcommand{\OS}{\mathcal {A}}
\def\ee{\emptyset}
\def\FF{{\mathcal F}}
\def\a{{\alpha}}
\def\b{{\beta}}
\def\d{{\delta}}

\newcommand{\s}{{\rm Sing\,}}

\newcommand{\End}{{\rm End\,}}
\newcommand{\Sing}{{\rm {Sing\,}}}
\newcommand{\Hess}{{\rm {Hess}^{(a)}}}

\newcommand{\V}{ V_{\bs \La}}

\newcommand{\Ll}{\La - \al(\bs k)}

\begin{document}

\title[Bethe Ansatz for Arrangements of Hyperplanes]
{Bethe Ansatz for Arrangements of Hyperplanes \\ and the Gaudin Model}

\author[{}]{
Alexander Varchenko ${}^{*,1}$}

\thanks{${}^1$ Supported in part by NSF grant DMS-0244579}

\begin{abstract}
We show that the Shapovalov norm of a Bethe vector in the Gaudin model
is equal to the Hessian of the logarithm of the corresponding master function
at the corresponding isolated critical point. We show that
different Bethe vectors are orthogonal. These facts are corollaries of
a general Bethe ansatz type construction, suggested in this paper and
 associated with an arbitrary arrangement of hyperplanes.

\end{abstract}

\maketitle 
\centerline{\it ${}^{*}$Department of Mathematics, University of
  North Carolina at Chapel Hill,} \centerline{\it Chapel Hill, NC
  27599-3250, USA} \medskip

\centerline{July, 2004}

\section{Introduction}

The Bethe ansatz is a large collection of methods in the theory of
quantum integrable models to calculate the spectrum and eigenvectors
for a certain commutative sub-algebra of observables for an integrable
model. Elements of the sub-algebra are called Hamiltonians, or
integrals of motion, or conservation laws of the model.  The
bibliography on the Bethe ansatz method is enormous, see for example
\cite{BIK, Fa, FT}.

In the theory of the Bethe ansatz one assigns the Bethe ansatz
equations to an integrable model. Then a solution of the Bethe ansatz
equations gives an eigenvector of commuting Hamiltonians of the
model. The general conjecture is that the constructed vectors form a
basis in the space of states of the model.

The simplest and interesting example is the Gaudin model associated
with a complex simple Lie algebra $\g$, see \cite{B, BF, F,
FFR, G, MV2, MV3, MV4, RV, ScV, V2, V3}.
One considers highest weight $\g$-modules
$V_{\La_1}, \dots , V_{\La_n}$ and their tensor product $\V$. One fixes
a point $z = (z_1, \dots , z_n) \in \C^n$ with distinct coordinates
and defines linear operators
$K_1(z), \dots , K_n(z)$ on $V_{\bs \La}$ by the formula
\bea
K_i(z)\ = \ \sum_{j \neq i}\ \frac{\Omega^{(i,j)}}{z_i - z_j} , 
\qquad i = 1, \dots , n .
\eea
Here $\Omega^{(i,j)}$ is the Casimir operator acting in the 
$i$-th and $j$-th factors of the tensor product. The operators are 
called  the Gaudin Hamiltonians
of the Gaudin model associated with $\V$. The Hamiltonians commute.

The common eigenvectors of the Gaudin Hamiltonians are constructed by
the Bethe ansatz method. Namely, one assigns to the model a scalar
function $\Phi (t,z)$ of new auxiliary variables $t$ and a $\V$-valued
function $\omega(t,z)$ such that $\omega(t^0,z)$ is an eigenvector of
the Hamiltonians if $t^0$ is a critical point of $\Phi$.  The
functions $\Phi$ and $\omega$ were introduced in \cite{SV} to
construct hypergeometric solutions of the KZ equations.  The function
$\Phi$ is called the { master function} and the function $\omega$ is
called { the canonical weight function}.

The first question is if the Bethe eigenvector $\omega(t^0,z)$ is
non-zero. In this paper we show that the Bethe vector is non-zero if
$t^0$ is a non-degenerate critical point of the master function
$\Phi$.  To show that we prove (in part (i) of Theorem \ref{main}) the
following identity: 
\bea\label{first} 
S ( \omega (t^0,z), \omega
(t^0,z) )\ =\ \text{Hess}_t\ \text{ln}\ \Phi (t^0, z) \ .  
\eea 
Here $S$ is the tensor product Shapovalov form on the tensor product $\V$ and the
right hand side of the formula is the Hessian at $t^0$ of the function
${\rm ln}\ \Phi$. 

This formula for the Gaudin model, associated with
$\g=sl_{r+1}$, was
proved in \cite{V2}, if $r = 1$, and for arbitrary $r$ in \cite{MV4},
see also \cite{Ko, R, RV, TV, MV1}.

We also show (in part (ii) of Theorem \ref{main}) that different Bethe
vectors are orthogonal with respect to the tensor product Shapovalov form.

These two statements allow us to reduce the Bethe ansatz conjecture to
a question about the number of non-degenerate
critical points of the master function, see part (iii)
of Theorem \ref{main}.

\bigskip

The formulated statements on the Bethe vectors are corollaries of a general 
construction, suggested in this paper and  
related to an arbitrary arrangement of hyperplanes. 
Namely, let $\A$ be an arrangement
of affine hyperplanes in $\C^k$ having a vertex. One defines the Orlik-Solomon algebra
$\OS(\A) = \oplus_p \OS^p(\A)$ and the flag space 
$\FF(\A) = \oplus_p \FF^p(\A)$ in the standard way, \cite{SV}.

The spaces $\OS^p(\A)$ and $\FF^p(\A)$ are dual. We are interested in the top degree
spaces  $\OS^k(\A)$ and $\FF^k(\A)$.

Assume that a complex number $a(H)$ is assigned to every hyperplane $H$ of $\A$. Then
one can define a symmetric bilinear form $S^{(a)} : \FF^k(\A)\otimes \FF^k(\A) \to \C$
called the Shapovalov form of $\A$, \cite{SV}. One also defines the master function
of $\A$, \ $\Phi = \prod_{H\in\A} f_H^{a(H)}$, where $f_H = 0$ is the defining equation
of the hyperplane $H$.

Let $t_1, \dots , t_k$ be coordinates in $\C^k$. Remind that the 
 space $\OS^k(\A)$ is the space
of rational differential $k$-forms on $\C^k$ which can be written as exterior
polynomials in differential 1-forms $df_H/f_H,\ H \in \A$. Hence each $\eta \in \OS^k(\A)$
can be written as $u \,dt_1\wedge\dots\wedge dt_k$
where $u$ is a rational function. 

Define the rational map $v : \C^k \to \FF^k(\A)$, regular on the complement to the
union of hyperplanes,  as follows. Let
$\epsilon \in \OS^k(\A) \otimes \FF^k(\A)$ be the canonical element,\
$ \epsilon = \sum_m x_m^*\otimes x_m$ where $\{x_m\}$ is a basis in
$\FF^k(\A)$ and $\{x^*_m\}$ is the dual basis in $\OS^k(\A)$.
If $x_m^* = u_m dt_1\wedge\dots\wedge dt_k$, then
$v(t) = \sum_m u_m(t)\, x_m$.

In part (ii) of Theorem \ref{first theorem} we show that
$$
S^{(a)}(v(t),v(t)) \ =\ (-1)^k\
\det_{1\leq i,j \leq k} 
( \frac{\partial^2} {\partial t_i \partial t_j} \ln \Phi ) (t) \ .
$$
In part (iii) of Theorem \ref{first theorem} we show that if
 $t^1, t^2$ are different isolated critical points of $\Phi$, 
then the special vectors $v(t^1), v(t^2)$ are orthogonal, \
$S^{(a)}(v(t^1),v(t^2))\ =\ 0$. 

Theorem \ref{first theorem} is the main result of the paper. 
To obtain the results concerning the Bethe ansatz for the Gaudin model we apply
Theorem \ref{first theorem} to discriminantal arrangements following methods of
\cite{SV}.

\bigskip

In this paper we considered the Bethe ansatz associated with a simple
Lie algebra.  In the same way one may consider the case of an
arbitrary Kac-Moody algebra. The statements and proofs remain the
same.

\bigskip

The paper is organized as follows.  Section 2 contains basic facts
about the Orlik-Solomon algebra and flag spaces of an
arrangement. Section 3  contains
the construction of special singular vectors in the top flag
space and the statement of Theorem \ref{first theorem}.
In Section 4 we prove  Theorem \ref{first theorem}.  Section 5
contains applications of Theorem \ref{first theorem} to the Bethe ansatz
associated with the Gaudin model.

\bigskip 

The idea of this paper was formulated long time ago in \cite{V2}, where it was
mentioned 
that an analog of the Bethe ansatz construction must exist for an arbitrary arrangement 
of hyperplanes.

\bigskip

The author thanks IHES for warm hospitality.

\section{Arrangements, \cite{SV, V1}}

\subsection{Arrangement}
Let $\A =\{H_j\}$, $j\in J(\A)$,  be an arrangement of affine hyperplanes in the
complex affine space $\C^k$. Denote by $U$  the complement to the union of all
hyperplanes,
$$
U = \C^k - \cup_{j\in J(\A)} H_j \ .
$$
Hyperplanes $H_j$ define in $\C^k$ the structure
of a stratified space.
A  closed  stratum $X_\a \subset \C^k$ is the intersection of some
hyperplanes  $H_j$, $j\in J_\a \subset J(\A)$. 
For a stratum $X_\a$ we denote
$l(X_\a) = \mathrm{codim}_{\C^k} X_\a$.
In this paper we will always assume that $\A$ has a vertex, a stratum of dimension 0.

\subsection{Orlik-Solomon algebra}
Define complex vector spaces $\OS^p(\A)$, $p = 0,  \dots, k$.
 For $p=0$ set $\OS^p(\A)=\C$. For  $p \geq 1$,\
 $\OS^p(\A)$   is generated by symbols
$(H_{i_1},...,H_{i_p})$ with $H_{i_j}\in\A$, such that
\begin{enumerate}
\item[(i)] $(H_{i_1},...,H_{i_p})=0$
if $H_{i_1}$,...,$H_{i_p}$ are not in general position, that is if the
intersection $H_{i_1}\cap ... \cap H_{i_p}$ is empty or
 its codimension is
 less than $p$;
\item[(ii)]
$ (H_{i_{\sigma(1)}},...,H_{i_{\sigma(p)}})=(-1)^{|\sigma|}
(H_{i_1},...,H_{i_p})
$
for any permutation $\sigma\in S_p$;
\item[(iii)]
$\sum_{j=1}^{p+1}(-1)^j (H_{i_1},...,\widehat{H}_{i_j},...,H_{i_{p+1}}) = 0
$
for any $(p+1)$-tuple $H_{i_1},...,H_{i_{p+1}}$ of hyperplanes
in $\A$ which are
not in general position and such that $H_{i_1}\cap...\cap H_{i_{p+1}}\not=\ee$.
\end{enumerate}

The direct sum $\OS(\A) = \oplus_{p=1}^{N}\OS^p(\A)$ is a
graded skew commutative algebra with respect to the
 multiplication
 $$
(H_{i_1},...,H_{i_p})\cdot(H_{i_{p+1}},...,H_{i_{p+q}}) =
 (H_{i_1},...,H_{i_p},H_{i_{p+1}},...,H_{i_{p+q}})\ .
$$
 The algebra is  called {\it the Orlik-Solomon algebra} of the arrangement $\A$.

Let $a: \A\to \C$ be a map which assigns to each hyperplane $H$
a complex number $a(H)$ called {\it the exponent} of $H$.
 Set
 $$
\omega(a)\ =\ \sum_{H\in\A}\,a(H)\,H\ {}\ \in \ {}\ \OS^1(\A)\ .
$$
 The multiplication by $\omega(a)$ defines a differential
 $$
d_\OS^{(a)}\ :\   \OS^p(\A)\ \to\ \OS^{p+1}(\A) ,
\qquad
 x \ \mapsto\ \omega(a)\cdot x \ ,
$$
in the  vector space of the
Orlik-Solomon algebra.

It is known that for generic exponents $a$,\ $H^p(\OS^\bullet(\A),d_\OS^{(a)}) = 0$ if $ p < k$
and \linebreak
dim $H^p(\OS^\bullet,d_\OS^{(a)}) = |\chi(U)|$, where $\chi (U)$ is the Euler 
characteristics of $U$, see \cite{A, STV}.

\subsection{Space of Flags }
 For a stratum $X_\a$, \ $l(X_\a)=p$, {\it a flag starting at $X_\a$} is a sequence
$$
X_{\a_0}\supset
X_{\a_1} \supset \dots \supset X_{\a_p} = X_\a
$$
of strata such that
$ l(X_{\a_j}) = j$ for $j = 0, \dots , p$.

{}  For a stratum $X_\a$,
 we define $\overline{\FF}_{X_\a}$  as  the
complex vector space  with basis vectors
$$
\overline{F}_{X_{\a_0},\dots,X_{\a_p}=X_\a}
$$
 la\-bel\-ed by the elements of
the set of all flags  starting at $X_\a$.

 Define  $\FF_{X_\a}$ as the quotient of
$\overline{\FF}_{X_\a}$ over the subspace generated by the vectors
\begin{equation}
\label{flagrelation} \sum\limits_{X_\b,\ 
X_{\a_{j-1}}\supset X_\b\supset X_{\a_{j+1}}}\ 
\overline {F}_{X_{\a_0},\dots,
X_{\a_{j-1}},X_{\b},X_{\a_{j+1}},\dots,X_{\a_p}=X_\a}\ .
\notag
\end{equation}
Such a vector is determined  by  $j \in \{ 1, \dots , p-1\}$ and
an incomplete flag $X_{\a_0}\supset...\supset
X_{\a_{j-1}} \supset X_{\a_{j+1}}\supset...\supset
X_{\a_p} = X_\a$ with $l(X_{\a_i})$ $=$ $i$.

Denote by ${F}_{X_{\a_0},\dots,X_{\a_p}}$ the image in $\FF_\a$ of the basis vector
$\overline{F}_{X_{\a_0},\dots,X_{\a_p}}$.  Set
$$
{\FF}^p(\A)\ =\ \oplus_{X_\a,\, l(X_\a)=p}\ {\FF}_{X_\a}\ ,
\qquad
{\FF}(\A)\ =\ \oplus _{p=0}^k\,{\FF}^p(\A)\ .
$$
The direct sum
 $$
{\FF}(\A)\ =\ \oplus_{p=0}^k\,{\FF}^p(\A)
$$
is a complex
with respect to the differential
 \bea
d_{\FF} :  \FF^p(\A) \to \FF^{p+1}(\A) ,
\qquad
{F}_{X_{\a_0},\dots,X_{\a_p}} \mapsto \sum_{X_{\a_{p+1}},\ 
X_{\a_{p}}\supset X_{\a_{p+1}}}
{F}_{X_{\a_0},\dots,X_{\a_p},X_{\a_{p+1}}}\ .
\eea

\subsection{Duality} 
The vector spaces $\OS^p(\A)$ and $\FF^p(\A)$ are dual. 
The pairing $ \OS^p(\A)\otimes\FF^p(\A) \to \C$ is defined as follows.
{}For $H_{i_1},...,H_{i_p}$ in general position, set
$F(H_{i_1},...,H_{i_p})=F_{X_{\a_0},\dots,X_{\a_p}}$
where
$$
X_{\a_0}=\C^k,\quad X_{\a_1}=H_{i_1},\quad \dots , \quad
X_{\a_p}=
H_{i_1} \cap \dots \cap H_{i_p} .
$$
Then set $\langle (H_{i_1},...,H_{i_p}), F \rangle = (-1)^{|\sigma|},$
if $F = F(H_{i_{\sigma(1)}},...,H_{i_{\sigma(p)}})$ for some $\sigma \in S_p$,
and $\langle (H_{i_1},...,H_{i_p}), F \rangle = 0$ otherwise.

Define the map $\delta_{\FF}^{(a)} :  \FF^p(\A)\to\FF^{p-1}(\A)$ to be
 the map adjoint to
$d_\OS^{(a)} :   \OS^{p-1}(\A) \to \OS^{p}(\A)$.

An element $v \in \FF^k(\A)$ will be called {\it  singular } if
$\delta_\FF^{(a)} v = 0$. Denote by $\Sing\,\FF^k(\A) \subset \FF^k(\A)$ the subspace
of singular vectors.

For generic exponents $a$ the dimension of $\Sing\FF^k(\A)$ is equal to $|\chi(U)|$.

\subsection{The Shapovalov map and form}
 The collection of exponents $a$ determines {\it the Shapovalov  map}
 $$
\mathcal S^{(a)} : \FF(\A) \to \OS(\A),
\quad
  {F}_{X_{\a_0},\dots,X_{\a_p}}\ \mapsto \
\sum \ a(H_{i_1}) \cdots a(H_{i_p})\ (H_{i_1}, \dots , H_{i_p})\ ,
$$
 where the sum is taken over all $p$-tuples $(H_{i_1},...,H_{i_p})$ such that
$$
H_{i_1} \supset X_{a_1},\ {}\ .\ .\ .\ {} , \ {} H_{i_p}\supset X_{\a_p}\ .
$$
According to \cite{SV}, the map $\mathcal S^{(a)}$ is a morphism of the
 complex $(\FF^\bullet(\A), d_{\FF})$ to the complex
$(\OS^\bullet(\A), d_{\OS}^{(a)})$. The image
$(\mathcal S^{(a)}(\FF^\bullet(\A)), d^{(a)}_\OS)$ 
is called {\it the complex of flag forms} of $\A$.

Identifying $\OS(\A)$ with $\FF(\A)^*$, we may consider
the map $\mathcal S^{(a)}$ as a bilinear form on the vector space $\FF(\A)$.
This bilinear form, denoted by $S^{(a)}$,
 is symmetric and  is called {\it the Shapovalov form}.

If $F_1, F_2 \in \FF^p(\A)$, then
\bea
S^{(a)}(F_1,F_2) = 
\sum_{\{i_1, \dots , i_p\} \subset J(\A)} \ a(H_{i_1}) \cdots a(H_{i_p})
\ \langle (H_{i_1}, \dots , H_{i_p}), F_1 \rangle
\ \langle (H_{i_1}, \dots , H_{i_p}), F_2 \rangle \ ,
\eea
where the sum is over all unordered $p$-element sets.

\section{Master Function, Special Vectors}

\subsection{Master function}\label{master function}
 Let $\A =\{H_j\}$, $j\in J(\A)$,  be an arrangement 
of affine hyperplanes in  $\C^k$.
For $j\in J(\A)$, fix a defining equation for $H_j$, $f_{j} = 0$.
Let $a: \A\to \C$ be a set of exponents. Then the function
$$
\Phi = \prod_{j\in J(\A)} f_j^{a(H_j)}
$$
is called {\it the master function}.
The master function is a multi-valued function defined on $U$.

A point $t\in U$ is called {\it a critical point} 
of $\Phi$ \ if\ $d \Phi\vert_t = 0$.

Fix affine coordinates $t_i, i = 1, \dots , k$, on $\C^k$. 
For $j\in J(\A)$, we have
$$
f_j(t_1, \dots , t_k) = b^0_j + b^1_j t_1 + \dots + b^k_j t_k ,
\qquad 
b^i_j  \in \C .
$$ 
A point $t\in U$ is  a critical point of $\Phi$ if
and only if
$$
\sum_{j\in J(\A)}
\ \frac {\partial f_j}{\partial t_i}\ 
\frac {a(H_j)}{f_j}\ = \ 0 \ , 
\qquad
i\ =\ 1, \dots k\ ,
$$
at $t$.

It is known that for generic exponents $a$ all critical points of
$\Phi$ are non-degenerate and their number is equal to $|\chi(U)|$,
see \cite{V2, OT, Si}.

\subsection{Realization of the Orlik-Solomon algebra}
For $j \in J(\A)$, consider the logarithmic differential form
$\omega_j = df_j/f_j$ on $\C^k$. 
Let $\bar{\OS}(\A)$ be the graded $\C$-algebra with unit element
generated by all $\omega_j$'s. 
The map ${\OS}(\A) \to \bar{\OS}(\A), \ H_j \mapsto \omega_j$, 
is an isomorphism. We shall identify ${\OS}(\A)$ and $\bar{\OS}(\A)$.

\subsection{Special vectors in $\FF^k(\A)$}\label{Special}

A top degree form $\eta \in {\OS}^k(\A)$ can be written as
$$
\eta \ = \ u \ dt_1 \wedge \dots \wedge dt_k
$$
where $u$ is a rational function  regular on $U$.

Define the rational map $v : \C^k \to \FF^k(\A)$, regular on U, as follows.
For $t\in U$,  set $v(t)$ to be the element of $\FF^k(\A)$ such that
$$
\langle \,\eta,\, v(t)\, \rangle\ =\ u(t)
\qquad
{\mathrm{for\ any\ }} \ \eta \in \OS^k(\A) \ .
$$
 Let
$\epsilon \in \OS^k(\A) \otimes \FF^k(\A)$ be the canonical element,\
$ \epsilon = \sum_m x_m^*\otimes x_m$ where $\{x_m\}$ is a basis in
$\FF^k(\A)$ and $\{x^*_m\}$ is the dual basis in $\OS^k(\A)$.
If $x_m^* = u_m dt_1\wedge\dots\wedge dt_k$, then
$v(t) = \sum_m u_m(t)\, x_m$.

The map $v$ will be called {\it the specialization map}, its value $v(t)$ 
will be called {\it the special vector} associated with $t\in U$.

Define the rational function $\Hess : \C^k \to \C$, regular on $U$, by the formula
$$
\Hess (t)\ =\ \det_{1\leq i,j \leq k} 
( \frac{\partial^2} {\partial t_i \partial t_j} \ln \Phi ) (t) \ .
$$

\begin{theorem} \label{first theorem}
${}$
\begin{enumerate}
\item[(i)] A point $t\in U$ is a critical point of $\Phi$, if and only if the special vector
$v(t)$ is a singular vector.
\item[(ii)]
If $t\in U$, then \
$$
S^{(a)}(v(t),v(t)) \ =\  (-1)^k\ \Hess (t)\ .
$$
\item[(ii)] If $t^1, t^2 \in U$ are different isolated critical points of $\Phi$, 
then the special singular vectors $v(t^1), v(t^2)$ are orthogonal, \
$$
S^{(a)}(v(t^1),v(t^2))\ =\ 0\ .
$$
\end{enumerate}
\end{theorem}

The theorem is proved in Section \ref{proofs}.

\subsection{Shapovalov images of special vectors}. 
The composition of the specialization 
and  Shapovalov maps define a rational map
$\mathcal S^{(a)}v : \C^k \to \OS^k(\A)$ regular on $U$.

\begin{lemma}\label{Sv lemma}
 For $t^1, t^2 \in U$, we have
\bea
\mathcal S^{(a)}v (t^1)|_{t^2}\ = \
S^{(a)}(v(t^1),v(t^2))\ dt_1\wedge \dots\wedge dt_k \ .
\eea
\end{lemma}
\begin{proof}
Let $t^1,t^2 \in U$ and
$\mathcal S^{(a)}v(t^1) \,=\, u\, dt_1\wedge\dots\wedge dt_k$. Then
$S^{(a)}(v(t^1), v(t^2)) = 
\langle \mathcal S^{(a)}v(t^1) , v(t^2) \rangle = u(t^2)$.
\end{proof}

The following corollary gives an estimate from above on the number of non-degenerate 
critical points of $\Phi$ by the dimension of the vector space
\bea
\mathcal{H}^k\ = \ \mathcal S^{(a)}(\FF^{k}(\A)) \,/\, d^{(a)} \mathcal
S^{(a)}(\FF^{k-1}(\A))\ .
\eea

\begin{corollary}\label{S span}
Let $\A$ be an arrangement of affine hyperplanes
in $\C^k$. Let $\rm{d}$ be a natural number. 
Let $C$ be a set of $\rm{d}$ non-degenerate critical points of $\Phi$.
Then the natural projection of vectors $\{ \mathcal S^{(a)}v(t) \}_{t \in C}$ 
 span a $\rm{d}$-dimensional subspace in $\mathcal H^k$.
\end{corollary}

In particular, if $\Phi$ has $\rm{d}$  non-degenerate
critical points, where $\rm{d}$ is the dimension
of $\mathcal H^k$, then  the natural projections to $\mathcal H^k$
of the vectors $ \mathcal S^{(a)}v(t)$, associated with those points, 
form a basis in $\mathcal H^k$.

\subsection{Basis of special singular vectors} 
The following corollary gives an estimate from above on the number of non-degenerate 
critical points of $\Phi$ by the dimension of the kernel of the linear operator
$\delta_{\FF}^{(a)} :  \FF^k(\A)\to\FF^{k-1}(\A)$.

\begin{corollary}\label{span}
Let $\A$ be an arrangement of affine hyperplanes
in $\C^k$. Let $\rm{d}$ be a natural number. 
Let $C$ be a set of $\rm{d}$ non-degenerate critical points of $\Phi$.
Then the special singular vectors $\{ v(t) \}_{t \in C}$  span a $\rm{d}$-dimensional
subspace in $\Sing \FF^k(\A)$.
\end{corollary}

In particular, if $\Phi$ has $\rm{d}$  non-degenerate
critical points, where $\rm{d}$ is the dimension
of $\Sing \FF^k(\A)$. Then  the special singular vectors, associated to those 
points, form  a basis in $\Sing \FF^k(\A)$.

\begin{corollary}
If the exponents $a$ are generic, then the set $\{ v(t) \}_{t \in C}$ 
is a basis in $\Sing \FF^k(\A)$.
\end{corollary}

\subsection{Arrangements with symmetries} Assume that a finite group $G$ acts
on $\C^k$ by affine linear transformations 
so that the arrangement $\A$ is preserved. Assume that exponents
$a$ are preserved by this action,  $a( g(H) )\,=\, a(H)$
for $g \in G$, $H\in \A$.

The group  $G$ naturally acts on $\FF^p(\A)$ for any $p$.
The action on $\FF^k(\A)$ 
will be denoted by  $R$. 
The action commutes with the differential $\d_\FF^{(a)}$.
The subspace $\Sing \FF^k(\A) \subset \FF^k(\A)$ is $G$-invariant. 
The Shapovalov form $S^{(a)} : \FF^k(\A) \otimes \FF^k(\A) \to \C$ 
is $G$-invariant.

Let $\Omega^k$ be the one dimensional complex vector space of differential $k$-forms on
$\C^k$ invariant with respect to all affine translations. 
The action of $G$ on $\C^k$ determines a representation 
$\rho : G \to \C^*, g \mapsto \rho_g$,
defined by the condition
$$
\rho_g\, g^*(\eta)\, = \,\eta\ ,
\qquad
 \eta \in \Omega^k\ .
$$ 
Fix affine coordinates $t_i, i = 1, \dots , k$, on $\C^k$. 
Let $v : \C^k \to \FF^k(\A)$ be the specialization map.
We have
\bea\label{action}
 v( g(t) )\, = \,\rho_g \, R_g ( v(t) ) \ ,
\qquad 
t\in U,\ g \in G \ .
\eea
The critical set $C \subset U$ of the master function $\Phi$ is $G$-invariant and
$$
\Hess (g(t))\ =\ (\rho_g)^2 \ \Hess (t) \ .
$$

\begin{corollary}

Let $t\in U$ be a non-degenerate critical point of $\Phi$ and
$\mathcal{O}$ its $G$-orbit.  Let $W$ be the span in $\Sing \FF^k(\A)$
of the vectors $\{ v(t') \}\vert_{t'\in \mathcal{O}}$. Then $W$ is
$G$-invariant and dim $W = |\mathcal{O}|$.

\end{corollary}

Let $\rho^1, \dots , \rho^N$ be all distinct irreducible
representations of $G$,\ $d_1 \dots , d_N$ the corresponding
dimensions, $\chi_1, \dots , \chi_N$ the corresponding characters,\
$\FF^k(\A) = W_1 \oplus \dots \oplus W_N$ the corresponding canonical
decomposition of $\FF^k(\A)$ into isotypical components.  The
projection $p_j$ of $\FF^k(\A)$ onto $W_j$ associated with this
decomposition is given by the formula \cite{S}
$$
p_j \ = \ \frac{d_j}{|G|} \ \sum_{g \in G}\ (\chi_j(g))^\dagger\  R_g \ ,
$$
where $z^\dagger$ denotes the complex conjugate of $z \in \C$.
Let $\OS^k(\A) = V_1 \oplus \dots \oplus V_N$
be the decomposition dual to $\FF^k(\A) = W_1 \oplus \dots \oplus W_N$.

For $j = 1, \dots , N$,  define the rational map $v_j : \C^k \to W_j$,
regular on $U$, as the composition of $v$ and $p_j$.

Let $\{x_m\}$ be a basis in $W_j$ and $\{x^*_m\}$  the dual basis in
$V_j$.  If $x_m^* = u_m dt_1\wedge\dots\wedge dt_k$, then $v_j(t) =
\sum_m u_m(t)\, x_m$. Clearly $v_1(t) + \dots + v_N(t) = v(t)$ for $t\in
U$.

The map $v_j$ will be called {\it the specialization map} associated with the isotypical
component $W_j \subset \FF^k(\A)$.


\begin{corollary}\label{main corollary}
${}$

\begin{enumerate}
\item[(i)]
Let $t^1, t^2 \in U$ be isolated critical points of $\Phi$ whose
$G$-orbits do not intersect. Then \
$S^{(a)} (v_j(t^1), v_j(t^2))\ =\ 0$ .

\item[(ii)]
Let $t\in U$ be an isolated critical point of $\Phi$. Assume that
the $G$-orbit of $t$ consists of $|G|$ elements. Then for $j = 1, \dots , N$, 
we have 
$$
S^{(a)} (v_j(t), v_j(t))\ =\ c_j \  (-1)^k\ \Hess (t) \ , 
\qquad
c_j\ = \ \frac{(d_j)^2}{|G|^2}\ \sum_{g\in G} 
((\chi_j(g))^\dagger)^2 \ .
$$
In particular, if $\chi_j :G \to \C$ takes values in $\R$ only, then
$c_j = (d_j)^2/|G|$.

\end{enumerate}
\end{corollary}

\section{Proof of Theorem \ref{first theorem}} \label{proofs}
\subsection{Proof of parts (i) and (ii) of Theorem \ref{first theorem}}
A point $t \in U$ is a critical point of $\Phi$ if and only if
the differential 1-form
$$
\omega^{(a)} \ = \ \sum_{j \in J(\A)}\ a(H_j)\ \frac{d f_j}{f_j}
$$
equals zero at $t$. The form is zero at $t$ if and only if
$\langle v(t) , \eta \rangle = 0$ for all $\eta$ lying in the image of
$d^{(a)}_\OS$, thus if and only if the vector $v(t)$ is singular.

Let $t^1,t^2 \in U$. By definition of the Shapovalov form, we have
\bea
S^{(a)}(v(t^1), v(t^2)) \ =\  
\sum_{\{j_1, \dots , j_k\} 
\subset J(\A)}\ D(j_1,\dots , j_k)^2\ 
\prod_{l=1}^k \
\frac{a(H_{j_l})}
{f_{j_l}(t^1)\ f_{j_l}(t^2)} \ ,
\eea
where \ $D(j_1,\dots , j_k) = \det_{1\leq i,l \leq k} (b^i_{j_l})$ and the
sum is over all unordered $k$-element subsets in $J(\A)$.
The right hand side of this 
formula for $t^1 = t^2 = t$  gives  $(-1)^k$ 
det $( \frac{\partial^2} {\partial t_i \partial t_j} \ln \Phi ) (t)$.

\subsection{Generic arrangements}
An arrangement $\A$ is generic if for any
distinct $i_1, \dots , i_k \in J(\A)$, the intersection
 $H_{i_1}\cap \dots \cap H_{i_k}$ is a point, and for any
distinct $i_1, \dots , i_{k+1} \in J(\A)$, the intersection
 $H_{i_1}\cap \dots \cap H_{i_{k+1}}$ is empty.

Fix an ordering on $J(\A)$.

For a generic arrangement, a basis in $\FF^k(\A)$ is formed by the
flags $F(H_{i_1},\dots,H_{i_k})$, \ $i_1, \dots , i_k \in J(\A)$,\
such that $i_1<\dots<i_k$. This basis will be called {\it standard}.

In $\FF^k(\A)$ we  have 
$$
F(H_{i_1},\dots,H_{i_k}) = (-1)^{|\sigma|}
F(H_{i_{\sigma(1)}},\dots,H_{i_{\sigma(k)}})
$$ 
for any  $\sigma \in S_k$.
 
We have
$$
S^{(a)}(F(H_{i_1},\dots,H_{i_k}) , F(H_{i_1},\dots,H_{i_k})) = a(H_{i_1})\cdots a(H_{i_k})
$$
and
$$
S^{(a)}(F(H_{i_1},\dots,H_{i_k}) , F(H_{j_1},\dots,H_{j_k})) = 0 
$$
for distinct elements of the standard basis.

For any distinct $j_1, \dots , j_{k+1} \in J(\A)$, $j_1 < \dots <
j_{k+1}$,\ define a linear map $L_{j_1, \dots , j_{k+1}} : \FF^k(\A)
\to \FF^k(\A)$ by its action on the elements of the standard basis: \
if $i_1, \dots , i_k$ is not a subset of $j_1, \dots , j_{k+1}$, then
$F(H_{i_1},\dots,H_{i_k}) \mapsto 0$, and
$$
F(H_{j_1}, \dots , \widehat{H_{j_p}}, \dots , H_{j_{k+1}}) \mapsto
(-1)^p \sum_{l=1}^{k+1} (-1)^l a(H_{j_l})
F(H_{j_1}, \dots , \widehat{H_{j_l}}, \dots , H_{j_{k+1}}) .
$$

\begin{lemma}
The map $L_{j_1, \dots , j_{k+1}}$ is self-adjoint, 
$$
S^{(a)}(L_{j_1, \dots , j_{k+1}} F_1 , F_2) =
S^{(a)}(F_1 , L_{j_1, \dots , j_{k+1}} F_2)
$$
for any $F_1, F_2 \in \FF^k(\A)$.
\hfill $\square$
\end{lemma}

Fix affine coordinates $t_i,\ i = 1, \dots , k$,\ on $\C^k$ and
for  $j\in J(\A)$  a polynomial
$$
f_j(t_1, \dots , t_k) =  b^0_j + b^1_j t_1 + \dots + b^k_j t_k 
$$ 
whose kernel is $H_j$.

Consider $f_j,\ j\in J(\A),$ 
as polynomials in variables $t_1, \dots , t_k, \ b^0_j,\ j\in J(\A)$.
For \linebreak
$j_1, \dots , j_{k+1} \in J(\A),\ j_1 < \dots < j_{k+1},$ introduce polynomials
$$
f_{j_1, \dots , j_{k+1}} \ =\ \sum_{p=1}^{k+1}\ D(j_1, \dots , \widehat{j_{p}},  
\dots , j_{k+1}) \ b^0_{j_p} \ .
$$ 
The polynomials $f_{j_1}, \dots , f_{j_{k+1}}, f_{j_1, \dots ,
j_{k+1}}$ are linearly dependent.  Denote $\omega_{j_1, \dots ,
j_{k+1}} = {d f_{j_1, \dots , j_{k+1}}}/ {f_{j_1, \dots , j_{k+1}}}$.
Then
$$
\omega_{j_1} \wedge \dots \wedge \omega_{j_{k+1}} =
\omega_{j_1, \dots , j_{k+1}} \wedge \sum_{p=1}^{k+1} (-1)^{p-1}
\omega_{j_1} \wedge \dots \wedge \widehat{\omega_{j_{p}}} \wedge \dots \wedge
\omega_{j_{k+1}} .
$$

\begin{lemma}\label{2}
\bea
&&
\sum_{j_1 < \dots < j_k }\ \big( \sum_{j} a(H_j)\
 \omega_j  \big) \wedge \omega_{j_1} \wedge \dots \wedge \omega_{j_k}
\otimes
F(H_{j_1}, \dots , H_{j_k}) =
\\
&&
\phantom{aaaa}
\sum_{j_1 < \dots < j_k } 
\sum_{i_1 < \dots < i_{k+1}}
\omega_{i_1, \dots , i_{k+1}} \wedge 
 \omega_{j_1} \wedge \dots \wedge \omega_{j_k}  
\otimes
L_{i_1, \dots , i_{k+1}} 
F(H_{j_1}, \dots , H_{j_k})\ .
\eea
\hfill
$\square$
\end{lemma}

For $j\in J(\A)$, define a linear map $K_j : \FF^k(\A) \to \FF^k(\A)$ by the formula
$$
K_j\ =\ 
\sum \ (-1)^p \frac{D(i_1, \dots , \widehat{i_{p}},  \dots , i_{k+1})}
{f_{i_1, \dots , i_{k+1}}} \ L_{i_1, \dots , i_{k+1}}
$$
where the sum is over all $i_1, \dots , i_{k+1} \in J(\A), \
i_1 < \dots < i_{k+1},$ and $1\leq p \leq k+1$ such that $i_p = j$.
The operator $K_j$ is self-adjoint.

\begin{lemma} If $t \in U$ is a critical point of $\Phi$, then for any 
$j \in J(\A)$, the special
singular vector $v(t)$ is an eigenvector of $K_j$ with eigenvalue
 $a(H_j)/f_j \vert_t$.
\end{lemma}

The lemma follows from Lemma \ref{2}

\begin{corollary}\label{orth corollary}
If $t^1, t^2 \in U$ are distinct critical points of $\Phi$, then
$v(t^1)$ and $v(t^2)$ are orthogonal with respect to $S^{(a)}$.

\end{corollary}

Now part (iii) of Theorem \ref{first theorem} follows from Corollary
\ref{orth corollary} and the continuity of \linebreak
$S^{(a)}(v(t^1),v(t^2))$ with respect to deformations of $t^1, t^2$
and of the arrangement $\A$.

\section{Applications to the Bethe ansatz of the Gaudin model}

\subsection{The Gaudin model}\label{model}
Let $\g$ be a simple Lie algebra over $\C$ with Cartan matrix
$A=(a_{i,j})_{i,j=1}^r$. 
Let $\h \subset \g$  be the Cartan sub-algebra.
Fix simple roots $\al_1, \dots , \al_r$ in $\h^*$
and an invariant bilinear form $( , )$ on $\g$. 
 Let $H_1, \dots , H_r \, \in \h$ be the corresponding
coroots, $\langle\la , H_i\rangle = 2 (\la,\al_i)/ (\al_i,\al_i)$
 for $\la\in\h^*$.
In particular, $\langle \al_j , H_i \rangle = a_{i,j}$.

Let $E_1, \dots , E_r\, \in \n_+,\ H_1, \dots , H_r\, \in \h,\
F_1, \dots , F_r\, \in \n_- $ be the Chevalley generators of $\g$,
\bea
[E_i, F_j] & = & \delta_{i,j} \,H_i ,
\qquad i, j = 1, \dots r ,
\\
{}[ h , h'] & = & 0 ,
\qquad 
\phantom{aaaaa}
h, h' \in \h ,
\\
{}[ h, E_i] &=& \langle \alpha_i, h \rangle\, E_i ,
\qquad  
h \in \h, \ i = 1, \dots r ,
\\
{}[ h, F_i] &=& - \langle \alpha_i, h \rangle\, F_i ,
\qquad  h\in \h, \ i = 1, \dots r ,
\eea
and \ $(\mathrm{ad}\,{} E_i)^{1-a_{i,j}}\,E_j = 0,$\ 
$(\mathrm{ad}\, {} F_i)^{1-a_{i,j}}\,F_j = 0,$ \
for all $i\neq j$.

Let $(x_i)_{i\in I}$ be an orthonormal basis in $\g$, \
$\Omega =  \sum_{i\in I} x_i\otimes x_i\ \in \g \otimes \g$
the Casimir element.

For a $\g$-module $V$ and $\mu \in \h^*$
denote by $V[\mu]$ the weight subspace of $V$ of weight $\mu$ and by
$\s V[\mu]$ the subspace of singular vectors of weight $\mu$,
\bea
\s V[\mu]\ =\ \{ \ v \in V\ |\ \n_+v = 0, \ hv = \langle \mu, h 
\rangle v \ \} \ .
\eea

Let $n$ be a positive integer and  $\bs \La = (\La_1, \dots , \La_n)$,
$\La_i \in \h^*$, a set of weights.
For $\mu \in \h^*$ let $V_{\mu}$ 
be the irreducible $\g$-module with highest weight $\mu$.
Denote by $V_{\bs \La}$ the tensor product 
$V_{\La_1} \otimes \dots \otimes V_{\La_n}$.

If $X \in \End\,(V_{\La_i})$, then we denote by $X^{(i)} \in \End
(V_{\bs \La})$ the operator $ \cdots \otimes \id \otimes X \otimes \id
\otimes \cdots$ acting non-trivially on the $i$-th factor of the
tensor product only. If $X = \sum_m X_m \otimes Y_m \in
\End (V_{\La_i} \otimes V_{\La_j})$, then we set 
$X^{(i,j)} = \sum_m X^{(i)}_m \otimes Y^{(j)}_m\ \in \End (\V)$.

Let $z = (z_1, \dots , z_n)$ be a point in $\C^n$ with distinct
coordinates.  Introduce linear operators $K_1(z), \dots , K_n(z)$ on
$V_{\bs \La}$ by the formula \bea K_i(z)\ = \ \sum_{j \neq i}\
\frac{\Omega^{(i,j)}}{z_i - z_j}\ , \qquad i = 1, \dots , n .  \eea
The operators are called {\it the Gaudin Hamiltonians} of the Gaudin
model associated with $\V$. The
Hamiltonians commute, $ [ K_i(z), K_j(z) ] = 0$ for all $i, j$.

The main problem for the Gaudin model is to diagonalize simultaneously
the  Hamiltonians.

One can check that the  Hamiltonians commute with the action of $\g$ on
$V_{\bs \La}$. Therefore it is enough to
diagonalize the  Hamiltonians on the subspaces of singular vectors
$\s V_{\bs \La}[\mu] \subset \V$. 

The eigenvectors of the Gaudin Hamiltonians are constructed by the
Bethe ansatz method.  We remind the construction in the next section.

\subsection{Master functions  and the canonical weight function, 
c.f. \cite{MV4}} \label{master sec} Fix a collection of weights
$\bs\La = (\La_1, \dots , \La_n)$, $\La_i\in \h^*$, and a collection
of non-negative integers $\bs k = (k_1, \dots , k_r)$.  Denote $k =
k_1 + \dots + k_r$, $\La = \La_1 + \dots + \La_n$, and $\al(\bs k) =
k_1\al_1 + \dots + k_r\al_r$.

Let $c$ be the unique non-decreasing function from $\{1, \ldots , k\}$
to $\{1, \ldots , r\}$, such that $\# c^{-1}( i) = k_i$ for $i = 1,
\dots , r$.  The {\it master function} $\Phi(t, z, \bs \La, \bs l)$
associated with this data is defined by the formula 
\bea\label{master}
\Phi(t, z, \bs \La, \bs k) = \prod_{1\leq i < j\leq n} (z_i -
z_j)^{(\La_i, \La_j)} \prod_{i=1}^l \prod_{s=k}^n (t_i -
z_s)^{-(\al_{c(i)}, \La_s)} \prod_{1 \leq i < j \leq k} (t_i -
t_j)^{(\al_{c(i)},\al_{c(j)})} , 
\eea 
see \cite{SV}.  The function
$\Phi$ is a function of complex variables $t = ( t_1, \dots , t_k)$,
$z = (z_1, \dots , z_n)$, weights $\bs \La$, and discrete parameters
$\bs k$.  The main variables are $t$, the other variables will be
considered as parameters.

For given $z,  \bs \La, \bs k$, a point $t \in \C^k$ is  {\it a critical
point} of the master function
if the following system of algebraic equations is satisfied,
\bea\label{Bethe eqn}
- \sum_{s=1}^n \frac{(\alpha_{c(i)}, \La_s)}{t_i - z_s}\ +\
\sum_{j,\ j\neq i} \frac{(\alpha_{c(i)}, \alpha_{c(j)})}{ t_i - t_j}
= 0, 
\qquad
i = 1, \dots , k .
\eea

Let $\Sigma_k$ be the permutation group of the set $\{1, \dots , k\}$.
Denote by $\bs \Sigma_{\bs k} \subset \Sigma_k$ the subgroup of all permutations preserving
the level sets of the function $c$. The subgroup $\bs\Sigma_{\bs k}$ is isomorphic to 
$\Sigma_{k_1}\times \dots \times \Sigma_{k_r}$ and 
acts on $\C^k$ permuting coordinates of $t$. The
action of the subgroup $\bs\Sigma_{\bs k}$ preserves the  critical set
of the master function. All orbits of the action of
$\bs\Sigma_{\bs k}$ on the critical set
have the same cardinality $k_1! \cdots k_r!$\ .

Consider highest weight irreducible $\g$-modules
$V_{\La_1}, \dots , V_{\La_n}$, the tensor product $\V
= V_{\La_1}\otimes\dots \otimes V_{\La_n}$, 
and its weight subspace $\V[\La - \al(\bs k)]$. Fix a highest weight vector
$v_{\La_i}$ in $V_{\La_i}$ for  all $i$.

We construct a rational map
\bea
\omega \ :\ \C^k \times \C^n\ \to \V[\Ll]
\eea
called {\it the canonical weight function}.

Let $P(\bs k,n)$ be the set of sequences 
$I\ = \ (i_1^1, \dots , i^1_{j_1};\ \dots ;\ i^n_1, \dots , i^n_{j_n})$ of integers in 
$\{1, \dots , r\}$ 
such that for all
$i = 1, \dots ,  r$, the integer $i$ appears in $I$ precisely $k_i$ times. 
For $I \in P(\bs k, n)$, and a permutation $\sigma \in \Sigma_k$, 
set $\sigma_1(i) = \sigma(i)$ for $i = 1, \dots , j_1$,
and $\sigma_s(i) = \sigma(j_1+\cdots +j_{s-1}+i)$ for $s = 2, \dots , n$ and 
$ i = 1, \dots , j_s$.

Define 
\bea
\Sigma(I)\ {} = \ {}
\{\ \sigma \in \Sigma_k\ {} |\ {} c(\sigma_s(j)) = i_s^j 
\ {} \text{for} \ {} s = 1, \dots , n \ {} \text{and} \ {} j = 1, \dots j_s\ \}\ .
\eea 

To every $I \in P(\bs k, n)$ we associate a vector 
\bea
F_Iv\ =\ F_{i_1^1} \dots F_{i_{j_1}^1}v_{\La_1} \otimes \cdots
\otimes F_{i_1^n} \dots F_{i_{j_n}^n}v_{\La_n}
\eea
in $\V[\Ll]$, and  rational functions
\bea
\omega_{I,\sigma} \ =\ \omega_{\sigma_1(1),\ldots,\sigma_1(j_1)}(z_1)\
\cdots\
\omega_{\sigma_n(1),\ldots,\sigma_n(j_n)}(z_n) ,
\eea
labeled by $\sigma\in \Sigma(I)$,
where 
\bea
\omega_{i_1,\ldots, i_j}(z_s) \ =\ \frac 1 {(t_{i_1}-t_{i_2}) \cdots
(t_{i_{j-1}}-t_{i_j}) (t_{i_j}-z_s)} .
\eea
We set
\bea\label{bethe vector}
\omega(z,t) \ =\  \sum_{I\in P(\bs k,n)}\ \sum_{\sigma\in \Sigma (I)}\ 
\omega_{I,\sigma}\ F_I v\ .
\eea

The canonical weight function  was introduced in \cite{SV}  to solve 
the KZ equations, see \cite{SV, FSV2, FMTV}. The hypergeometric solutions to the KZ equations
with values in $\s  \V[\Ll]$ have the form
\bea
I(z)\ =\ \int_{\gamma(z)} \ \Phi(t,  z, \bs \La, \bs k)^{1/\kappa} \ \omega(t,z)\ dt .
\eea
Different formulas for the canonical weight function see in \cite{RSV}.

The values of the canonical weight function at the critical points (with respect to 
variables $t$)
of the master function are called {\it the Bethe vectors}, 
see \cite{RV, V2, FFR}.

\begin{theorem}[\cite{RV}]\label{cr bethe}
Assume that $z \in \C^n$ has distinct coordinates.
Assume that  $t \in \C^k$ is a critical point of the master function
$\Phi(\, .\, ,  z, \bs \La, \bs k)$. Then 
the vector $\omega(t,z)$ belongs to $\s\V[\Ll]$ and
 is an eigenvector of the Gaudin Hamiltonians $K_1(z), \dots , K_n(z)$.
\end{theorem}

This theorem was proved in \cite{RV} using the  
quasi-classical asymptotics of the hypergeometric solutions
of the KZ equations. The theorem also follows 
directly from Theorem  6.16.2 in \cite{SV},
cf. Theorem 7.2.5 in \cite{SV},\ see also Theorem 4.2.2 in \cite{FSV2}.

\subsection{The Shapovalov Form}
Define the anti-involution $\tau : \g \to \g $ sending
$E_1, \dots , E_r, \linebreak H_1, \dots , H_r, \
F_1, \dots , F_r$ to
$F_1, \dots , F_r, \ H_1, \dots , H_r, \
E_1, \dots , E_r$, respectively.

Let $W$ be a highest weight $\g$-module   with highest weight vector $w$.
{\it The Shapovalov form}  on $W$ is the unique
symmetric bilinear form $S$
defined by the conditions:
\bea
S(w, w) = 1 ,
\qquad
S(xu, v) = S(u, \tau(x)v)
\eea
for all $u,v \in W$ and $x \in \g$, see \cite{K}.

Let $V_{\La_1}, \dots , V_{\La_n}$ be irreducible
highest weight modules and
$\V$ their tensor product.
Let $v_{\La_i} \in V_{\La_i}$ be 
a highest weight vector and 
$S_i$ the corresponding Shapovalov form on $V_{\La_i}$. 
 Define 
a symmetric bilinear form on $\V$  by the formula
\bean\label{shap}
S\ =\ S_1 \otimes \cdots \otimes S_n .
\eean
The form $S$ will be called {\it the tensor product Shapovalov form on $\V$}.

\subsection{Application of Theorem \ref{first theorem}}
As in Section \ref{master sec} fix 
 a collection of weights $\bs\La = (\La_1, \dots , \La_n)$, $\La_i\in  \h^*$,
and a collection of non-negative integers  $\bs k = (k_1, \dots , k_r)$.

Let $S_{\V}$ be the tensor product Shapovalov from on the tensor product $\V$.
 
Fix a collection of distinct complex numbers $z = (z_1, \dots , z_n)$.

Let $t^1$, $t^2 \in \C^k$ be points such that $t^1$ has distinct coordinates and
$t^2$ has distinct coordinates and such that
none of coordinates of $t^1, t^2$  belongs to the set $\{z_1, \dots , z_n\}$. 

Under these assumptions,  we have the following theorem.

\begin{theorem}\label{main}
${}$

\begin{enumerate}
\item[(i)] Assume that
$t^1$ and $t^2$ are isolated critical points of 
$\Phi(\, .\, ,  z, \bs \La, \bs k)$. 
Assume that the $\bs\Sigma_{\bs k}$-orbits of $t^1$ and $t^2$ do not 
intersect. 
 Then the Bethe vectors $\omega(z,t^1)$ and $\omega(z,t^2)$ are orthogonal with respect 
to the tensor product Shapovalov form, \ 
$S_{\V}( \omega(z,t^1) , \omega(z,t^2) ) = 0$.

\item[(ii)] Assume that $t^1$ is an isolated critical point of
$\Phi(\, .\, ,  z, \bs \La, \bs k)$. Then
$$
S_{\V}( \omega(z,t^1) , \omega(z,t^1) ) =  
\det_{1\leq i,j \leq k} 
( \frac{\partial^2} {\partial t_i \partial t_j} \ln \Phi (t^1, z, \bs \La, \bs k) ) \ .
$$
\item[(iii)]  Let $\rm{d}$ be a natural number. 
Let $C$ be a set of $\rm{d}$ distinct $\bs\Sigma_{\bs k}$-orbits 
of non-degenerate critical points of $\Phi$.
Choose a representative $t^i$ in each orbit. Assume that
each point $t^i$ has distinct coordinates and none of the coordinates of $t^i$ belongs to
the set $\{z_1, \dots , z_n\}$. Then the Bethe vectors 
\ $\omega(z,t^i)$,  $i = 1, \dots \rm{d}$, span a $\rm{d}$-dimensional
subspace  in  $\s\V[\Ll]$.
\end{enumerate}
\end{theorem}

Part (i) of the theorem
was proved in \cite{RV} for $\g=sl_2$. Part (ii) was proved for $\g=sl_{r+1}$, $r = 1$,
 in \cite{V2} and for arbitrary $r$ in \cite{MV4}. In all 
those cases the proof used asymptotics of
Bethe vectors in the asymptotic zone $|z_1 - z_2| << \dots << |z_1 - z_n|$.

Part (iii) gives a  bound from above on 
the number $d$ of orbits of non-degenerate critical points of 
$\Phi (\,.\,, z, \bs \La, \bs k)$ in terms of the representation theory.
In particular, if the weight $\Ll$ is not integral dominant, then
$\Phi (\,.\,, z, \bs \La, \bs k)$ does not have at all non-degenerate critical points
(since in that case the space  $\s\V[\Ll]$ has dimension zero).

It is interesting to note that if $\Ll$ is not integral dominant, then all critical
points of $\Phi (\,.\,, z, \bs \La, \bs k)$ are non-isolated and the connected components
of the critical set are isomorphic to suitable Bruhat cells of the flag variety of the
Langlands dual Lie algebra, see \cite{ScV, MV2, MV3}.

In \cite{ScV} the case of $\g=sl_2$ was considered. It was proved that
if the weight $\Ll$ is integral dominant and $z_1, \dots , z_n$ are generic, 
then the function $\Phi (\,.\,, z, \bs \La, \bs k) $ has non-degenerate critical points
only and the critical points form $d$ orbits, where $d$ is the dimension of $\s\V[\Ll]$.
In particular, this means that the corresponding Bethe vectors form a basis
in $\s\V[\Ll]$.

\subsection{Proof of Theorem \ref{main}} 
For given $z = (z_1, \dots , z_n)$, the 
 discriminantal arrangement $\A(z)$ in $\C^k$ is defined as the collection of 
hyperplanes
\bea
H_{i}^s\ :\ t_i - z_s = 0 \ (i = 1, \dots , k, \ s = 1, \dots , n )\ ,
\qquad
H_{i,j}\ :\ t_i - t_j = 0 \ (1 \leq i < j \leq k ) \ ,
\eea
see \cite{SV}.
Define the weights of $\A(z)$ by the rule,\ $a(H_i^s) = - (\al_i, \La_s)$,\
$a(H_{i,j}) = - (\al_i, \al_j)$. \ Then the master function 
$\Phi (\,.\,, z, \bs \La, \bs k) $ is equal up to a constant factor to
 the master function $\Phi$,
defined in Section \ref{master function} and 
assigned to the weighted arrangement $\A(z)$.

Let $S^{(a)} : \FF^k(\A(z))\otimes \FF^k(\A(z)) \to \C$ be the Shapovalov form of $\A(z)$.

The group $\bs \Sigma_{\bs k}$ acts on $\C^k$ permuting coordinates. The action preserves
the discriminantal arrangement and its weights. Hence the group acts on $\FF^k(\A(z))$
and $\OS^k(\A(z))$. Set
\bea
W^- = \{ x \in \FF^k(\A(z))\ | \ R_\sigma (x) = (-1)^{|\sigma|} x, \ \sigma \in 
\bs \Sigma_{\bs k} \}\ .
\eea
Similarly define $V^- \subset \OS^k(\A(z))$ to be the skew-symmetric part of
$\OS^k(\A(z))$.  The subspaces $V^-$ and $W^-$ are dual.

For an element 
$I\ = \ (i_1^1, \dots , i^1_{j_1};\ \dots ;\ i^n_1, \dots , i^n_{j_n})$ in
$P(\bs k,n)$ and a permutation $\sigma$ in $\Sigma(I)$ define a flag 
$ f_{I,\sigma} \in \FF^k(\A(z))$,
\bea
 f_{I,\sigma} = F( H^1_{\sigma_1(1)},\dots ,H^1_{\sigma_1(j_1)},\dots ,
H^n_{\sigma_n(1)},\dots ,H^n_{\sigma_1(j_n)}) \ ,
\eea
and then an element $f_{I}  \in W^-$,
\bea
f_{I} \ =\  \frac 1{k_1!\,\dots\,k_r!}\
\sum_{\sigma \in \Sigma(I)}\ (-1)^{|\sigma|}\ f_{I,\sigma} \ .
\eea
\begin{theorem}[ Theorem 6.6 in \cite{SV} ] \label{Shap}
For $I,J \in P(\bs k,n)$,
$$
S_{\V} ( F_I v , F_J v )\ =\ (-1)^k\  k_1!\,\dots\,k_n!\ S^{(a)} ( f_I  , f_J  ) \ .
$$
\end{theorem}
\begin{theorem} \label{canonical}
The element
$$
\sum_{I\in P(\bs k,n)}\  (\sum_{\sigma\in \Sigma (I)}\ 
\omega_{I,\sigma} dt_1\wedge \dots dt_k )\ \otimes \ F_I v
$$
is the canonical element in $V^-\otimes W^-$.
\end{theorem}
The theorem is a direct corollary of Theorems 5.13 and 6.16.2 in \cite{SV}.

Let $v^- : \C^k \to W^-$ be the specialization map associated with the isotypical 
component $W^-\subset \FF^k(\A(z))$.
Let $t^1, t^2$ be as in Theorem \ref{main}, then
\bea
S_{\V}( \omega(z,t^1) , \omega(z,t^2) )\ =\ (-1)^k\ k_1!\,\dots\,k_n!\
S^{(a)} ( v^-(t^1) ,  v^-(t^2) ) 
\eea
by Theorems \ref{Shap} and \ref{canonical}. 
By Corollary \ref{main corollary} the right hand side is zero if the orbits 
of $t^1$ and $t^2$ do not intersect. 
By Corollary \ref{main corollary}, we have
\bea
(-1)^k\ k_1!\,\dots\,k_n!\ S^{(a)} ( v^-(t^1) ,  v^-(t^1) ) \ =\ 
\det_{1\leq i,j \leq k} 
( \frac{\partial^2} {\partial t_i \partial t_j} \ln \Phi (t^1, z, \bs \La, \bs k) ) \ .
\eea
This proves parts (i) and (ii) of Theorem \ref{main}. Part (iii) follows from the fact
that vectors
\ $\omega(z,t^i)$,  $i = 1, \dots, d$, have non-zero Shapovalov norm and are pair-wise 
orthogonal.

\end{document}